\documentclass[12pt]{amsart}

\usepackage{latexsym}
\usepackage{amsthm,amsfonts,amssymb,amscd,mathrsfs}
\usepackage[leqno]{amsmath}
\usepackage{xspace}
\usepackage{cite}
\usepackage{fancyhdr}
\input{xypic}
\xyoption{all}

\newcommand{\numberseries}{\mdseries}	

\newlength{\thmtopspace}		
\newlength{\thmbotspace}		
\newlength{\thmheadspace}		
\newlength{\thmindent}			

\newtheoremstyle{bfupright head,slanted body}
		{\thmtopspace}{\thmbotspace}
		{\slshape}{\thmindent}{\bfseries}{.}{\thmheadspace}
		{{\numberseries \thmnumber{(#2) }}\thmnote{#3}}

\newtheoremstyle{bfupright head,upright body}
		{\thmtopspace}{\thmbotspace}
		{\upshape}{\thmindent}{\bfseries}{.}{\thmheadspace}
		{{\numberseries \thmnumber{(#2) }}\thmnote{#3}}

\newtheoremstyle{bfit head,upright body}
		{\thmtopspace}{\thmbotspace}
		{\upshape}{\thmindent}{\upshape}{.}{\thmheadspace}
		{{\numberseries\thmnumber{(#2) }}
		{\bfseries\itshape\thmnote{\negthickspace#3}}}

\newtheoremstyle{it head,upright body}
		{\thmtopspace}{\thmbotspace}
		{\upshape}{\thmindent}{\upshape}{.}{\thmheadspace}
		{{\numberseries\thmnumber{(#2) }}
		{\itshape\thmnote{\negthickspace#3}}}

\newtheoremstyle{fixed bf head,slanted body}
		{\thmtopspace}{\thmbotspace}{\slshape}
		{\thmindent}{\bfseries}{.}{\thmheadspace}
		{{\numberseries \thmnumber{(#2) }}\thmname{#1}\thmnote{ (#3)}}

\newtheoremstyle{fixed bf head,upright body}
		{\thmtopspace}{\thmbotspace}{\upshape}
		{\thmindent}{\bfseries}{.}{\thmheadspace}
		{{\numberseries \thmnumber{(#2) }}\thmname{#1}\thmnote{ (#3)}}

\newtheoremstyle{indented paragraph}
		{\thmtopspace}{\thmbotspace}
		{\upshape}{\thmindent}{\upshape}{}{0pt}
		{\thmnote{#3 }}

\theoremstyle{bfupright head,slanted body}
\newtheorem{res}{}[section]		\newtheorem*{res*}{}

\theoremstyle{bfit head,upright body}
			\newtheorem*{com*}{}

\theoremstyle{bfupright head,upright body}
\newtheorem{bfhpg}[res]{}		\newtheorem*{bfhpg*}{}

\theoremstyle{it head,upright body}
		\newtheorem*{ithpg*}{}

\theoremstyle{fixed bf head,slanted body}
\newtheorem{thm}[res]{Theorem}		\newtheorem*{thm*}{Theorem}
\newtheorem{prp}[res]{Proposition}	\newtheorem*{prp*}{Proposition}
	\newtheorem*{cor*}{Corollary}
\newtheorem{lem}[res]{Lemma}		\newtheorem*{lem*}{Lemma}

\theoremstyle{fixed bf head,upright body}
\newtheorem{dfn}[res]{Definition}	\newtheorem*{dfn*}{Definition}
	\newtheorem*{obs*}{Observation}
\newtheorem{rmk}[res]{Remark}		\newtheorem*{rmk*}{Remark}
\newtheorem{exa}[res]{Example}		\newtheorem*{exa*}{Example}
		\newtheorem*{exe*}{Exercise}
		\newtheorem{stp*}{Setup}
	\newtheorem*{dfns*}{Definitions}
	\newtheorem*{obss*}{Observations}
		\newtheorem*{rmks*}{Remarks}
	\newtheorem*{exas*}{Examples}

\theoremstyle{indented paragraph}

\newlength{\thmlistleft}	
\newlength{\thmlistright}	
\newlength{\thmlistpartopsep}	
\newlength{\thmlisttopsep}	
\newlength{\thmlistparsep}	
\newlength{\thmlistitemsep}	

\newcounter{eqc} 
\newenvironment{eqc}{\begin{list}{\upshape (\textit{\roman{eqc}})}%
		    {\usecounter{eqc}%
		        \setlength{\leftmargin}{\thmlistleft}%
			\setlength{\labelwidth}{\thmlistleft}%
			\setlength{\rightmargin}{\thmlistright}%
			\setlength{\partopsep}{\thmlistpartopsep}%
			\setlength{\topsep}{\thmlisttopsep}%
			\setlength{\parsep}{\thmlistparsep}%
			\setlength{\itemsep}{\thmlistitemsep}}}%
		    {\end{list}}%

\newcounter{prt}
\newenvironment{prt}{\begin{list}{\upshape (\alph{prt})}%
		    {\usecounter{prt}%
		        \setlength{\leftmargin}{\thmlistleft}%
			\setlength{\labelwidth}{\thmlistleft}%
			\setlength{\rightmargin}{\thmlistright}%
			\setlength{\partopsep}{\thmlistpartopsep}%
			\setlength{\topsep}{\thmlisttopsep}%
			\setlength{\parsep}{\thmlistparsep}%
			\setlength{\itemsep}{\thmlistitemsep}}}%
		    {\end{list}}%

\newcounter{rqm}
\newenvironment{rqm}{\begin{list}{\upshape (\arabic{rqm})}%
		    {\usecounter{rqm}%
		        \setlength{\leftmargin}{\thmlistleft}%
			\setlength{\labelwidth}{\thmlistleft}%
			\setlength{\rightmargin}{\thmlistright}%
			\setlength{\partopsep}{\thmlistpartopsep}%
			\setlength{\topsep}{\thmlisttopsep}%
			\setlength{\parsep}{\thmlistparsep}%
			\setlength{\itemsep}{\thmlistitemsep}}}%
		    {\end{list}}%

			{\end{list}}%

\newenvironment{itemlist}{\nopagebreak \begin{list}{$\bullet$}%
		       {\setlength{\leftmargin}{\thmlistleft}%
			\setlength{\labelwidth}{\thmlistleft}%
			\setlength{\rightmargin}{\thmlistright}%
			\setlength{\partopsep}{\thmlistpartopsep}%
			\setlength{\topsep}{\thmlisttopsep}%
			\setlength{\parsep}{\thmlistparsep}%
			\setlength{\itemsep}{\thmlistitemsep}}}%
			{\end{list}}%

			{\end{list}}%
\newlength{\myindent}
{\setlength{\myindent}{\parindent}\begin{list}{}%
			{\setlength{\leftmargin}{#1}\setlength{\rightmargin}{#1}%
			\setlength{\partopsep}{0pt}%
			\setlength{\topsep}{\thmtopspace}%
			\setlength{\parsep}{0pt}%
			\setlength{\itemsep}{0pt}}
			\item[]}
			{\end{list}}%

\newenvironment{proof*}{\begin{proof}}{\renewcommand{\qed}{} \end{proof}}

  \newcommand{\proofof}[2][:]{``#2''#1}

\newcommand{\dispand}[1][and]{\hbox to \hsize{#1 \hfill} \nonumber \\}

\newlength{\seqsplit}
\title{Relative Ext groups, resolutions, and Schanuel classes}

\author{Henrik Holm }

\address{ {\flushleft Department} of Mathematical Sciences, University
  of Aarhus, Ny Mun\-kegade, Building 530,
  DK--8000 Aarhus C, Denmark}
\email{holm@imf.au.dk}
\urladdr{http://home.imf.au.dk/holm/}

 

\newcommand{\Ker}{\operatorname{Ker}}
\newcommand{\F}{\mathsf{F}}
\renewcommand{\S}{\mathcal{S}}
\newcommand{\Hom}{\operatorname{Hom}}
\newcommand{\Ext}{\underline{\operatorname{Ext}}}
\newcommand{\ext}{\operatorname{Ext}}

\begin{document}

\begin{abstract} 
  Given a precovering (also called contravariantly finite) class $\F$
  there are three natural approaches to a homological dimension with
  respect to $\F$: One based on Ext functors relative to $\F$, one
  based on $\F$--resolutions, and one based on Schanuel classes
  relative to $\F$. In general these approaches do not give the same
  result.  In this paper we study relations between the three
  approaches above, and we give necessary and sufficient conditions
  for them to agree.
\end{abstract}

\maketitle
\setcounter{section}{-1}


\section{Introduction}

\noindent
The fact that the category of modules over any ring $R$ has enough
projectives is a cornerstone in classical homolo\-gi\-cal algebra. The
existence of enough projective modules has three important
consequences:

\begin{itemlist}
\item To every module $A$, and integer $n$ one can define the
  \emph{Ext functor},
  \begin{displaymath}
    \ext_R^n(-,A),
  \end{displaymath}
  with well-known properties, see \cite[chap.\,V]{CE}.

\item Every module $M$ admits a \emph{projective resolution},
  cf.~\cite[chap.\,V]{CE}: 
  \begin{displaymath}
    \cdots \longrightarrow P_2 \longrightarrow P_1 \longrightarrow P_0
    \longrightarrow M \longrightarrow 0.
  \end{displaymath}

\item Every module $M$ represents a projective equivalence class
  $[M]$, and to this one can associate its \emph{Schanuel class},
  \begin{displaymath}
    \mathcal{S}([M]) := [\Ker \pi],
  \end{displaymath}
  where \mbox{$\pi \colon P \longrightarrow M$} is any epimorphism and
  $P$ is projective. One can also consider the iterated Schanuel maps
  $\mathcal{S}^n(-)$ for $n \geqslant 0$, see Schanuel's lemma
  \cite[chap.\,4, thm.\,A]{IK}.
\end{itemlist}
The three fundamental types of objects described above --- Ext
functors, projective resolutions, and Schanuel classes --- are linked
together as nicely as one could hope for, in the sense of the
following well-known result (see \cite[chap.\,V, prop.\,2.1]{CE}):

\begin{res*}[Theorem A]
  For any $R$--module $M$, and any integer \mbox{$n
  \geqslant 0$} the following conditions are equivalent:
  \begin{eqc}
  \item $\ext_R^{n+1}(M,A)=0$ for all $R$--modules $A$.
  \item There exists a projective resolution for $M$ of length $n$,
    \begin{displaymath}
      0 \longrightarrow P_n \longrightarrow \cdots \longrightarrow P_0
      \longrightarrow M \longrightarrow 0.
  \end{displaymath}
  \item $\mathcal{S}^n([M])=[0]$.
  \end{eqc}
\end{res*}

\noindent
The equivalent conditions of the theorem above define what it means
for a module $M$ to have projective dimension $\leqslant n$.

\medskip
\noindent
In relative homological algebra one substitutes the class of
projective modules by any other \emph{precovering} class $\F$, see
\eqref{bfhpg:precoves}. The fact that $\F$ is precovering allows for
well-defined constructions of:
\begin{itemlist}
\item Ext functors $\Ext_\F^n(-,A)$ relative to $\F$;

\item $\F$--resolutions, $\cdots \longrightarrow F_2 \longrightarrow
  F_1 \longrightarrow F_0 \longrightarrow 0$; and

\item Schanuel maps $\mathcal{S}_\F^n(-)$ relative to $\F$.
\end{itemlist}
The constructions of these relative objects are well-known, see for
example \cite[chap.\,8]{EnochsRHA} and \cite[lem.\,2.2]{EEEOMGJLO},
but for the benefit of the reader we give a short recapitulation in
Section \ref{sec:Preliminaries}.

\medskip
\noindent
Now, one could hope that there might exist an ``$\F$--version'' of
Theorem A, indeed, one would need such a theorem to have a rich and
flexible notion of an $\F$--dimension. Unfortunately, Theorem A fails
for a general precovering class $\F$!  The aim of this paper is to
understand, for a given precovering class $\F$, the different kind of
obstructions which keep Theorem A from being true.


\medskip 
\noindent 
In Section \ref{sec:RvsE} we investigate how the Ext condition $(i)$
and the resolution condition $(ii)$ in the $\F$--version of Theorem A
are related:

\medskip
\noindent
It is trivial that \mbox{$(ii) \Rightarrow (i)$} holds always, so we
restrict our attention to the converse implication. In Lemma
\eqref{lem:EtoR} we give a necessary condition for \mbox{$(i)
  \Rightarrow (ii)$}. In Theorem \eqref{thm:EtoR} we give a sufficient
condition for \mbox{$(i) \Rightarrow (ii)$} in terms of almost epi
precovers.  We also give concrete examples of precovering classes for
which the implication \mbox{$(i) \Rightarrow (ii)$} fails, and others
for which it holds.

\medskip
\noindent
In Section \ref{sec:RvsS} we study how the resolution condition $(ii)$
and the Schanuel condition $(iii)$ in the $\F$--version of Theorem A
are related:

\medskip
\noindent
The main results are Theorems \eqref{thm:StoR} and \eqref{thm:RtoS}
which give necessary and sufficient conditions for the implication
\mbox{$(iii) \Rightarrow (ii)$}, respectively, \mbox{$(ii) \Rightarrow
  (iii)$}, to hold.  We also present concrete examples of precovering
classes for which the implications \mbox{$(iii) \Rightarrow (ii)$} and
\mbox{$(ii) \Rightarrow (iii)$} fail.

\section*{Acknowledgements}

\noindent
It is a pleasure to thank Edgar E. Enochs and Christian U. Jensen for
their useful comments.

\section{Preliminaries} \label{sec:Preliminaries}

\begin{bfhpg}[Setup] 
  \label{bfhpg:setup} 
  Throughout, $R$ will be a ring, and all modules will be left
  $R$--modules. We write $\mathsf{Mod}\,R$ for the category of (left)
  $R$--modules, and $\mathsf{Ab}$ for the category of abelian groups.

  \medskip
  \noindent
  $\F$ will be any precovering class of modules,
  cf.~\eqref{bfhpg:precoves} below, which contains $0$ and is closed
  under isomorphism and finite direct sums.
\end{bfhpg}

\begin{bfhpg}[Precovering classes] \label{bfhpg:precoves} For
  definitions and results on precovering classes we generally follow
  \cite[chap.\,5 and 8]{EnochsRHA}. We mention here just a few notions
  which will be important for this paper.

  \medskip
  \noindent
  Let $\F$ be a class of modules. An \emph{$\F$--precover} of a module
  $M$ is a homomorphism $F \longrightarrow M$ with $F \in \F$, such
  that given any other homomorphism $F' \longrightarrow M$ with $F'
  \in \F$ then there exists a factorization,
  \begin{displaymath}
    \xymatrix{
      {} & F' \ar[d] \ar@{-->}[dl] \\
      F \ar[r] & M.
    }
  \end{displaymath}
  If every module admits an $\F$--precover then $\F$ is called
  \emph{precovering}. An (augmented) \emph{$\F$--resolution} of a
  module $M$ is a complex (which is not necessarily exact),
  \begin{displaymath}
    \cdots \longrightarrow F_2 \stackrel{\partial_2}{\longrightarrow}
    F_1 \stackrel{\partial_1}{\longrightarrow} F_0 
    \stackrel{\partial_0}{\longrightarrow} M \longrightarrow 0,
  \end{displaymath}
  with $F_0,F_1,F_2,\ldots \in \F$, such that
  \begin{displaymath}
    \cdots \longrightarrow (F,F_2) \stackrel{(F,\partial_2)}{\longrightarrow}
    (F,F_1) \stackrel{(F,\partial_1)}{\longrightarrow} (F,F_0) 
    \stackrel{(F,\partial_0)}{\longrightarrow} (F,M) \longrightarrow 0
  \end{displaymath}
  is exact for all \mbox{$F \in \F$}. 
  When $\F$ is precovering, and $T \colon \mathsf{Mod}\,R
  \longrightarrow \mathsf{Ab}$ is a contravariant additive functor,
  then one can well-define the $n$'th right derived functor of $T$
  relative to $\F$,
  \begin{displaymath}
    \textnormal{R}^n_\F T \colon \mathsf{Mod}\,R \longrightarrow
    \mathsf{Ab}.
  \end{displaymath}
  One computes $\textnormal{R}^n_\F T(M)$ by taking an non-augmented
  $\F$--resolution of $M$, applying $T$ to it, and then taking the
  $n$'th cohomology group of the resulting complex. For a module $A$
  we write:
  \begin{displaymath}
    \Ext^n_\F(-,A) = \textnormal{R}^n_\F\Hom_R(-,A).
  \end{displaymath}
  Note that we underline the Ext for good reasons: There is also a
  notion of a \emph{preenveloping} class. If $\mathsf{G}$ is
  preenveloping then one can right derive the Hom functor in the
  covariant variable with respect to $\mathsf{G}$. Thus for each
  $R$--module $B$ there are functors
  $\overline{\ext}^{\,n}_{\mathsf{G}}(B,-)$. However, in general,
  \begin{displaymath}
    \Ext^n_\F(B,A) \neq \overline{\ext}^{\,n}_{\mathsf{G}}(B,A)
  \end{displaymath}
  even if \mbox{$\F=\mathsf{G}$} is both precovering and
  preenveloping.
\end{bfhpg}

\begin{bfhpg}[$\F$--equivalence]
  Two modules $K$ and $K'$ are called \emph{$\F$--equivalent}, and we
  write \mbox{$K \equiv_\F K'$}, if there exist $F, F' \in \F$ with $K
  \oplus F' \cong K' \oplus F$. We use $[K]$ to denote the
  $\F$--equivalence class containing $K$.

  \medskip
  \noindent
  Now let $M$ be any module. By the version of Schanuel's lemma found
  in \cite[lem.\,2.2]{EEEOMGJLO}, the kernels of any two
  $\F$--precovers of $M$ are $\F$--equivalent. Thus the class $[\Ker
  \varphi]$, where $\varphi \colon F \longrightarrow M$ is any
  $\F$--precover of $M$, is a well-defined object depending only on
  $M$. We write
  \begin{displaymath}
    \S_\F(M)=[\Ker \varphi].
  \end{displaymath}
  As $\F$ is closed under finite direct sums; cf.~Setup
  \eqref{bfhpg:setup}, it is not hard to see that $\S_\F(M)$ only
  depends on the $\F$--equivalence class of $M$, and hence we get the
  induced Schanuel map:
  \begin{displaymath}
    \xymatrix{\mathsf{Mod}\,R/\!\equiv_\F \ar[rr]^-{\S_\F} & &
      \mathsf{Mod}\,R/\!\equiv_\F}\!,
  \end{displaymath}
  For \mbox{$n>0$} we write $\S^n_\F$ for the $n$--fold composition of
  $\S_\F$ with itself, and we set $\S^0_\F = \mathrm{id}$.
\end{bfhpg}

\noindent
This paper is all about studying relations between the conditions from
the following definition.

\begin{dfn} \label{dfn:conditions} 
  For any module $M$ and any integer $n \geqslant 0$ we consider the
  conditions:
  \begin{eqc}
  \item[$(E_{M,n})$] $\Ext^{n+1}_\F(M,A)=0$ for all modules $A$.
  \item[$(R_{M,n})$] There exists an $\F$--resolution of the form
    \begin{displaymath}
      0 \longrightarrow F_n \longrightarrow \cdots \longrightarrow F_0
      \longrightarrow M \longrightarrow 0. 
    \end{displaymath}
  \item[$(S_{M,n})$] $\S^n_\F([M])=[0]$.
  \end{eqc}
\end{dfn}

\begin{rmk}
  The conditions in Definition \eqref{dfn:conditions} are labeled
  according to the following mnemonic rules: ``$E$'' is for Ext,
  ``$R$'' is for Resolution, and ``$S$'' is for Schanuel.
\end{rmk}

\section{Relative Ext functors and resolutions} \label{sec:RvsE}

\noindent
In this section we study how the Ext condition and the resolution
condition of Definition \eqref{dfn:conditions} are related. It is
straightforward, cf.~Proposition \eqref{prp:strightforward} below,
that the resolution condition implies the Ext condition. The converse
is, in general, not true, but in Theorem \eqref{thm:EtoR} we give a
sufficient condition on $\F$ for this to happen.

\begin{prp} \label{prp:strightforward} 
  For any precovering class $\F$ we have:
  \begin{displaymath}
    (R_{M,n}) \Longrightarrow (E_{M,n}) \text{ for all modules
      $M$ and all integers $n \geqslant 0$.} \qed
  \end{displaymath}
\end{prp}

\begin{exa} \label{exa:nodirectsummands} 
  There exist precovering classes which are not closed under direct
  summands:

  \medskip
  \noindent
  Let $R$ be a left noetherian ring which is not Quasi--Frobenius, and
  set \mbox{$D=R \oplus E$} where $E$ is any non-zero injective
  $R$--module.  Define $\F$ to be the class of all modules which are
  isomorphic to $D^{(\Lambda)}$ for some index set $\Lambda$ (here
  $D^\varnothing =0$).  Note that $\F$ is precovering as for example
  an $\F$--precover of a module $M$ is given by the natural map
  \begin{displaymath}
    D^{(\Hom_R(D,M))} \longrightarrow M.
  \end{displaymath}
  To see that $\F$ is not closed under direct summands we note that
  $E$ is a direct summand of \mbox{$D \in \F$}. However, there exists
  no set $\Lambda$ for which \mbox{$E \cong D^{(\Lambda)}$} (since $R$
  is a direct summand of $D^{(\Lambda)}$ for any $\Lambda \neq
  \varnothing$, and since $R$ is not self-injective).
\end{exa}

\noindent
The example above makes the following lemma relevant:

\begin{lem} \label{lem:EtoR}
  A necessary condition for $\F$ to satisfy the implication: 
  \begin{displaymath}
    \text{$(E_{M,0}) \Longrightarrow (R_{M,0})$ for all modules $M$,}
  \end{displaymath}
  is that $\F$ is closed under direct summands.
\end{lem}

\begin{proof}
  Assume that $\F$ is not closed under direct summands. Then there
  exists an \mbox{$F \in \F$} and a direct summand $M$ of $F$ with
  \mbox{$M \notin \F$}.  We claim that $(E_{M,0})$ holds but that
  $(R_{M,0})$ does not:
 
  \medskip
  \noindent
  As $M$ is a direct summand of $F$, and as $\F$ is closed under
  finite direct sums, cf.~Setup \eqref{bfhpg:setup}, the abelian group
  $\Ext^1_\F(M,A)$ is a direct summand of $\Ext^1_\F(F,A)$ for every
  module $A$. The latter is zero as $F \in \F$, and hence also
  $\Ext^1_{\F}(M,A)=0$.  Now suppose for contradiction that there do
  exist an $\F$--resolution of $M$ of length zero:
  \begin{displaymath}
    0 \longrightarrow F_0
    \stackrel{\partial_0}{\longrightarrow} M \longrightarrow
    0. 
  \end{displaymath}
  We claim that $\partial_0$ must be an isomorphism (contradicting the
  fact that $M \notin \F$). As $M$ is a direct summand of $F$ there is
  a canonical embedding \mbox{$\iota \colon M \longrightarrow F$} and
  a canonical projection \mbox{$\pi \colon F \longrightarrow M$} with
  $\pi\iota=\mathrm{id}_M$.  As $\partial_0$ is an $\F$--precover of
  $M$, we get a factorization:
  \begin{displaymath}
    \xymatrix{
      {} & F \ar@{->>}[d]^-{\pi} \ar@{-->}[dl]_-{\varphi} \\
      F_0 \ar[r]_-{\partial_0} & M
    }
  \end{displaymath}
  It follows that $\partial_0(\varphi\iota)=\pi\iota = \textrm{id}_M$,
  so $\partial_0$ is epi and the sequence
  \begin{align*}
    \tag{\text{$\dagger$}}
    \xymatrix{
      0 \ar[r] & \Ker \partial_0  \ar[r] & F_0 
      \ar@<0.5ex>[r]^-{\partial_0} 
      & M \ar[r] \ar@<0.5ex>@{-->}[l]^-{\varphi\iota} & 0
    }
  \end{align*}
  splits. By assumption, $\Hom_R(G,\partial_0)$ is mono for all
  \mbox{$G \in \F$}, so by $(\dagger)$ it follows that
  \mbox{$\Hom_R(G,\Ker
    \partial_0)=0$} for all \mbox{$G \in \F$}. In parti\-cu\-lar,
  \begin{displaymath}
    \Hom_R(F_0,\Ker \partial_0)=0,
  \end{displaymath}
  and therefore $\Ker
  \partial_0 =0$ since $\Ker \partial_0$ is a direct summand of $F_0$.
  Consequently, $\partial_0$ is an isomorphism.
\end{proof}

\begin{lem} \label{lem:eqc}
  For a homomorphism $\varphi \colon F \longrightarrow M$ the
  following two conditions are equivalent:
  \begin{prt}
  \item Every endomorphism $g \colon M \longrightarrow M$ with
    $g\varphi = \varphi$ is an automorphism.
  \item Every endomorphism $g \colon M \longrightarrow M$ with
    $g\varphi = \varphi$ admits a left inverse.
  \end{prt}
\end{lem}

\begin{proof}
  We only need to show that (b) implies (a). Thus assume that
  $g\varphi=\varphi$. By assumption (b) there exists a homomorphism $v
  \colon M \longrightarrow M$ with $vg=\mathrm{id}_M$. Now
  \begin{displaymath}
    v\varphi = vg\varphi = \mathrm{id}_M\varphi = \varphi,
  \end{displaymath}
  so another application of (b) gives that also $v$ has a left
  inverse. As $v$ has $g$ as a right inverse, $v$ must be an
  automorphisms with $v^{-1}=g$.
\end{proof}

\begin{dfn} \label{dfn:almostepi} 
  A homomorphism $\varphi \colon F \longrightarrow M$ satisfying the
  equi\-va\-lent conditions of Lemma \eqref{lem:eqc} is called
  \emph{almost epi}.
  The precovering class $\F$ is called \emph{precovering by almost
    epimorphisms} if every module has an $\F$--precover which is
  almost epi.
\end{dfn}

\begin{exa} \label{exa:domain} 
  Clearly, every epimorphism is almost epi, but the converse is in
  general not true, as for example 
  \begin{displaymath}
    \mathbb{Z} \stackrel{2\cdot}{\longrightarrow} \mathbb{Z}
  \end{displaymath}
  is an almost epimorphism of abelian groups. It follows from Lemma
  \eqref{lem:stability} below that if a precovering class contains all
  free modules, then it is precovering by almost epimorphisms.
\end{exa}

\begin{lem} \label{lem:stability}
  If there exists an almost epi homomorphism \mbox{$\varphi \colon F
    \longrightarrow M$} with \mbox{$F \in \F$} then every
  $\F$--precover of $M$ is almost epi.
\end{lem}

\begin{proof}
  If \mbox{$\tilde{\varphi} \colon \tilde{F} \longrightarrow M$} is
  any $\F$--precover of $M$ then there exists a factorization,
  \begin{displaymath}
    \xymatrix{
      {} & F \ar[d]^-{\varphi} \ar@{-->}[dl]_-{\psi} \\
      \tilde{F} \ar[r]_-{\tilde{\varphi}} & M.
    }
  \end{displaymath}
  For any endomorphism \mbox{$g \colon M \longrightarrow M$} with
  $g\tilde{\varphi} = \tilde{\varphi}$ it follows that
  \begin{displaymath}
    g\varphi = g\tilde{\varphi}\psi = \tilde{\varphi}\psi = \varphi,
  \end{displaymath}
  and hence $g$ must be an automorphism since $\varphi$ is almost epi.
\end{proof}

\noindent
The next result gives much more information than Example
\eqref{exa:domain}, namely that there do indeed exist module classes
$\F$ which are precovering by almost epimorphisms, without every
$\F$--precover being epi. We postpone the proof of Proposition
\eqref{prp:representaton-finite} to the end of this section.

\begin{prp} \label{prp:representaton-finite}
  Consider the local ring \mbox{$R=\mathbb{Z}/4\mathbb{Z}$}. We denote
  the generator \mbox{$2+4\mathbb{Z}$} of the maximal ideal by $\xi$,
  and the residue class field $R/(\xi) \cong \mathbb{F}_2$ by $k$.

  \medskip
  \noindent
  Furthermore, let \mbox{$\F = \mathsf{Add}\,k$} be the class of all
  direct summands of set-indexed coproducts of copies of $k$. Then the
  following hold:
  \begin{prt}
  \item $\F$ is precovering by almost epimorphisms, cf.~Definition
    \eqref{dfn:almostepi}.
  \item $R$ does not admit an epi $\F$--precover.
  \end{prt}
\end{prp}

\noindent
The reason we are interested in classes which are precovering by
almost epimorphisms is because of the next result:

\begin{thm} \label{thm:EtoR}
  Assume that $\,\F$ is closed under direct summands and is
  precovering by almost epimorphisms. Then 
  \begin{displaymath}
    \text{$(E_{M,n}) \Longrightarrow (R_{M,n})$ for all modules $M$ and
      all integers $n \geqslant 0$.}
  \end{displaymath}
\end{thm}

\begin{proof}
  First we deal with the case $n=0$: Thus let $M$ be any module, and
  assume that $\Ext_{\F}^1(M,A)=0$ for all modules $A$. We must prove
  the existence of an $\F$--resolution of $M$ of length zero,
  \begin{displaymath}
    0 \longrightarrow G_0 \longrightarrow M \longrightarrow 0. 
  \end{displaymath}  
  By assumption on $\F$ we can build an $\F$--resolution of $M$ by
  successively taking almost epi $\F$--precovers $\varphi_0, \varphi_1,
  \varphi_2, \ldots$:
  \begin{displaymath}
    \xymatrix@R=3ex@C=3ex{
      {} & 0 \ar[dr] & {} & 0 & {} & {} & {} & 0 & {} & {}
    \\ 
      {} & {} & K_1 \ar[dr]^-{i_1} \ar@{-->}[ur] & {} & {} & {} & M
      \ar@{-->}[ur] \ar@{=}[dr] &
      {} & {} & {} 
    \\
      \cdots \ar[r] & F_2 \ar[rr]^-{\partial_2} \ar[ur]^-{\varphi_2} &
      {} & F_1 \ar[rr]^-{\partial_1} \ar[dr]_-{\varphi_1} & {} & F_0
      \ar[rr]^-{\partial_0} \ar[ur]^-{\varphi_0}
      & {} & M \ar[rr] & {} & 0 
    \\
      {} & {} & {} & {} & K_0 \ar@{-->}[dr] \ar[ur]_-{i_0} & {} & {} &
      {} & {} & {} 
    \\
      {} & {} & {} & 0 \ar[ur] & {} & 0 & {} & {} & {} & {}
    }
  \end{displaymath}
  We keep in mind that the $\F$--precovers $\varphi_n$ are not
  necessarily epi, and this is the reason why some of the arrows in
  the diagram above have been dotted. Applying $\Hom_R(-,A)$, for any
  module $A$, to the $\Hom_R(\F,-)$ exact complex,
  \begin{displaymath}
    0 \longrightarrow K_0 \longrightarrow F_0 \longrightarrow M
    \longrightarrow 0, 
  \end{displaymath}
  induces by \cite[thm.\,8.2.3(2)]{EnochsRHA} an exact sequence of
  relative Ext groups,
  \begin{displaymath}
    \tag{\text{$*$}}
    \Ext_{\F}^0(F_0,A) \stackrel{q}{\longrightarrow} \Ext_{\F}^0(K_0,A)
    \longrightarrow \Ext_{\F}^1(M,A)=0. 
  \end{displaymath}
  As $F_0 \in \F$ we have $\Ext_{\F}^0(F_0,A)=\Hom_R(F_0,A)$. Furthermore,
  \begin{align*}
    \Ext_{\F}^0(K_0,A) &= \Ker \Hom_R(\partial_2,A) \\
    &= \{ f \in \Hom_R(F_1,A) \,|\, f \partial_2=0 \},
  \end{align*}
  and one verifies that the homomorphism $q$ is given by \mbox{$g
    \longmapsto g\partial_1$} for \mbox{$g \in \Hom_R(F_0,A)$}.
  Applying these considerations to \mbox{$A=K_0$} and considering
  \mbox{$\varphi_1 \in \Ext_{\F}^0(K_0,K_0)$}, exactness of $(*)$
  implies the existence of a $g \in \Hom_R(F_0,K_0)$ with
  $g\partial_1=\varphi_1$, that is, $g i_0 \varphi_1 = \varphi_1$. As
  $\varphi_1$ is almost epi, $g i_0 \colon K_0 \longrightarrow K_0$
  must be an automorphism, and hence the sequence
  \begin{align*}
    \xymatrix{
      0 \ar[r] & K_0 \ar@<0.5ex>[r]^-{i_0} & F_0 \ar[r]^-{\pi_0}
      \ar@<0.5ex>@{-->}[l]^-{(g i_0)^{-1}g} 
      & F_0/K_0 \ar[r] & 0
    }
  \end{align*}
  is split exact. In particular, \mbox{$F_0/K_0 \in \F$} as $\F$ is
  closed under direct summands, and we claim that the induced
  monomorphism,
  \begin{align*}
    \xymatrix@R=3ex@C=3ex{
      F_0/K_0 \ar@{^(-->}[rr]^-{\overline{\varphi}_0} & & M \\
      & F_0 \ar[ur]_-{\varphi_0} \ar@{->>}[ul]^-{\pi_0} & 
    }
  \end{align*}
  is an $\F$--precover of $M$. To see this let \mbox{$\varphi \colon F
    \longrightarrow M$} be a homomorphism with \mbox{$F \in \F$}. As
  \mbox{$\varphi_0 \colon F_0 \longrightarrow M$} is an $\F$--precover
  there exists \mbox{$\psi \colon F \longrightarrow F_0$} with
  $\varphi_0\psi = \varphi$. Consequently, \mbox{$\pi_0\psi \colon F
    \longrightarrow F_0/K_0$} satisfies
  $\overline{\varphi}_0(\pi_0\psi)=\varphi$.
  \begin{displaymath}
    \xymatrix@R=5ex@C=5ex{
      F_0 \ar[d]_-{\pi_0} \ar[dr]^-{\varphi_0} & F \ar[d]^-{\varphi}
      \ar@{-->}[l]_-{\psi} 
      \\ 
      F_0/K_0 \ar[r]_-{\overline{\varphi}_0} & M
    }  
  \end{displaymath}
  Now, as $\overline{\varphi}_0$ is a mono $\F$--precover of $M$,
  \begin{displaymath}
    0 \longrightarrow F_0/K_0
    \stackrel{\overline{\varphi}_0}{\longrightarrow} M \longrightarrow
    0 
  \end{displaymath}
  is an $\F$--resolution of $M$ of length zero.

  \medskip
  \noindent
  Finally we consider the case $n>0$: We assume that
  $\Ext_{\F}^{n+1}(M,A)=0$ for all modules $A$, and we must prove the
  existence of an $\F$--resolution of $M$ of length $n$.  Let
  \mbox{$\partial_0 \colon F_0 \longrightarrow M$} be an
  $\F$--precover of $M$. By \cite[thm.\,8.2.3(2)]{EnochsRHA} the
  $\Hom_R(\F,-)$ exact complex
  \begin{displaymath}
    \tag{\text{$\dagger$}}
    0 \longrightarrow \Ker \partial_0 \longrightarrow F_0
    \stackrel{\partial_0}{\longrightarrow} M 
    \longrightarrow 0
  \end{displaymath}  
  induces, for any module $A$, a long exact sequence of relative Ext
  groups:
  \begin{displaymath}
    0=\Ext_\F^n(F_0,A) \longrightarrow \Ext_\F^n(\Ker \partial_0,A)
    \longrightarrow \Ext_\F^{n+1}(M,A)=0.
  \end{displaymath}
  It follows that 
  \begin{displaymath}
    \Ext_\F^{(n-1)+1}(\Ker \partial_0,A) = \Ext_\F^n(\Ker
    \partial_0,A) = 0,  
  \end{displaymath}
  so the induction hypothesis implies that $\Ker \partial_0$ admits an
  $\F$--resolution of length $n-1$, say,
  \begin{displaymath}
   \tag{\text{$\ddagger$}}
    0 \longrightarrow F_n \longrightarrow \cdots \longrightarrow F_1
    \longrightarrow  \Ker \partial_0 \longrightarrow 0.  
  \end{displaymath}  
  Pasting together $(\dagger)$ and $(\ddagger)$ we get the desired
  $\F$--resolution of $M$ of length $n$.
\end{proof}

\begin{proof}[Proof of Proposition \eqref{prp:representaton-finite}]
  Note that $R$ is a two-dimensional $k$--vector space with basis
  $\{1,\xi\}$, so every element of $R$ has a unique representation of
  the form $a+b\xi$ where $a,b \in k \cong \mathbb{F}_2$.

  \medskip
  \noindent
  Just as in Example \eqref{exa:nodirectsummands} it follows that $\F
  = \mathsf{Add}\,k$ is precovering, but we shall also prove this more
  directly below.

  \medskip
  \noindent
  It is useful to observe that a homomorphism \mbox{$F \longrightarrow
    M$} with \mbox{$F \in \F$} is an $\F$--precover of $M$ if and only
  if every homomorphism \mbox{$k \longrightarrow M$} admits a
  factorization:
  \begin{displaymath}
    \tag{\text{$\natural$}}
    \xymatrix{ {} & k \ar@{-->}[dl] \ar[d] \\
     F \ar[r] & M.
    }
  \end{displaymath}
  One important consequence of this is that if \mbox{$F_j
    \longrightarrow M_j$} is a family of $\F$--precovers then the
  coproduct \mbox{$\coprod_j F_j \longrightarrow \coprod_j M_j$} is
  again an $\F$--precover.

  \medskip
  \noindent
  For every $c \in k$ there is an $R$--linear map
  \begin{displaymath}
    \varphi_c \colon k \longrightarrow R \, , \ a \longmapsto ac\xi,
  \end{displaymath}
  and it is not hard to see that, in fact, every $R$--linear map
  \mbox{$k \longrightarrow R$} has the form $\varphi_c$ for some $c
  \in k$.  Combining this with the commutative diagram
  \begin{displaymath}
    \xymatrix{ {} & k \ar[dl]_-{c\cdot} \ar[d]^-{\varphi_c} \\
     k \ar[r]_-{\varphi_1} & R,
    }
  \end{displaymath}
  observation $(\natural)$ implies that \mbox{$\varphi_1 \colon k
    \longrightarrow R$} is an $\F$--precover of $R$. Since $\varphi_1$
  is not epi, $R$ cannot be the homomorphic image of any module from
  $\F$, and this proves (b) from the proposition.

  \medskip
  \noindent
  We are now ready to prove part (a) of the proposition, namely that
  every $R$--module admits an almost epi $\F$--precover. It is
  well-known\footnote{The author is convinced that this result and its
    natural generalizations must be folklore.  However, since the
    author was not able to find a reference, a quick argument is given
    below.

    \medskip
    \noindent
    Let $M \neq 0$ be any $R$--module. Since
    \mbox{$R=\mathbb{Z}/4\mathbb{Z}$} only has the two proper ideals
    $(0)$ and $(\xi)$ there are two possibilities:
    \begin{rqm}
    \item for every $0 \neq x \in M$ we have
      $\mathrm{Ann}_R(x)=(\xi)$, or
    \item there exists $0 \neq x \in M$ with $\mathrm{Ann}_R(x)=(0)$.
    \end{rqm}
    In case (1) we can consider $M$ as a module over the field
    $k=R/(\xi)$, and it follows that $M \cong k^{(I)}$ for some index
    set $I$. In case (2) there is a monomorphism $R \longrightarrow
    M$, and since $R$ is self-injective it follows that $R$ is a
    direct summand of $M$. Using Zorn's lemma we find a maximal free
    (=injective) direct summand $R^{(J)}$ of $M$, and hence we can
    write $M=M' \oplus R^{(J)}$ where $M'$ satisfies condition (1).}
  that every $R$--module is isomorphic to one of the form
  \begin{displaymath}
    k^{(I)} \oplus R^{(J)}
  \end{displaymath}
  for suitable index sets $I$ and $J$.  Hence we only need to show that the module $k^{(I)}
  \oplus R^{(J)}$ has an almost epi $\F$--precover. By the observation
  $(\natural)$ it follows that
  \begin{displaymath}
    \xymatrix@C=20ex@R=-0.5ex{k^{(I)} & k^{(I)} \\
     \oplus \ar[r]^-{\varphi = \scriptsize
     \left(\!\!
       \begin{array}{ll}
         \mathrm{id}_{k^{(I)}} & 0 \\
         0 & \varphi_1^{(J)}
       \end{array}
     \!\!\!\right)
     } & \oplus \\
     k^{(J)} & R^{(J)}
    }
  \end{displaymath}
  is an $F$--precover. To argue that $\varphi$ is almost epi we let 
  \begin{displaymath}
    \xymatrix@C=20ex@R=-0.5ex{k^{(I)} & k^{(I)} \\
     \oplus \ar[r]^-{g= \scriptsize
     \left(\!\!
       \begin{array}{ll}
         g_{11} & g_{12} \\
         g_{21} & g_{22}
       \end{array}
     \!\!\!\right)
     } & \oplus \\
     R^{(J)} & R^{(J)}
    }
  \end{displaymath}
  be any endomorphism with \mbox{$\varphi = g\varphi$}. We must prove
  that $g$ is an automorphism. By assumption,
  \begin{align*}
     \tag{\text{$*$}}
     { \small \left(\!\!
       \begin{array}{ll}
         \mathrm{id}_{k^{(I)}} & 0 \\
         0 & \varphi_1^{(J)}
       \end{array}
     \!\!\!\right) }
     =
     \left(\!\!
       \begin{array}{ll}
         g_{11} & g_{12} \\
         g_{21} & g_{22}
       \end{array}
     \!\!\!\right)  
     {\small \left(\!\!
       \begin{array}{ll}
         \mathrm{id}_{k^{(I)}} & 0 \\
         0 & \varphi_1^{(J)}
       \end{array}
     \!\!\!\right) }
     =
     { \left(\!\!
       \begin{array}{ll}
         g_{11} & g_{12}\,\varphi_1^{(J)} \\
         g_{21} & g_{22}\,\varphi_1^{(J)}
       \end{array}
     \!\!\!\right) }
  \end{align*}
  In particular it follows that \mbox{$g_{11}=\mathrm{id}_{k^{(I)}}$}
  and \mbox{$g_{21}=0$}, so $g$ takes the form
  \begin{align*}
     g =
     \left(\!\!
       \begin{array}{ll}
         \mathrm{id}_{k^{(I)}} & g_{12} \\
         0 & g_{22}
       \end{array}
     \!\!\!\right). 
  \end{align*}
  If we can prove that $g_{22} \colon R^{(J)} \longrightarrow R^{(J)}$
  is an automorphism, then $g$ must be an automorphism as well with
  inverse
  \begin{align*}
     g^{-1} =
     \left(\!\!
       \begin{array}{ll}
         \mathrm{id}_{k^{(I)}} & -g_{12}\,g_{22}^{-1} \\
         0 & g_{22}^{-1}
       \end{array}
     \!\!\right).  
  \end{align*}
  To see that $g_{22}$ is an automorphism we use another relation from
  $(*)$, namely that \mbox{$\varphi_1^{(J)}=g_{22}\,\varphi_1^{(J)}$}.
  As
  \begin{displaymath}
    g_{22} \in \Hom_R(R^{(J)},R^{(J)}) \cong \left( R^{(J)} \right)^J 
  \end{displaymath}
  it follows --- if we consider the elements of $R^{(J)}$ as
  $J$--columns --- that $g_{22}$ is given by multiplication from the
  left by a unique $J\times J$--matrix $(r_{ij})_{i,j \in J}$ with
  entries from $R$, in which each column $(r_{ij})_{i\in J}$ belongs
  to $R^{(J)}$. More precisely, $g_{22}$ is given by the formula:
  \begin{displaymath}
    R^{(J)} \ni \{s_j\}_{j\in J} \longmapsto (r_{ij})_{i,j \in
      J}\{s_j\}_{j\in J} = \left\{ \sum_{j\in J}
      r_{ij} s_j \right\}_{i \in I} \in R^{(J)}.
  \end{displaymath}
  Of course, \mbox{$\varphi_1^{(J)} \colon k^{(J)} \longrightarrow
    R^{(J)}$} is given by the $J\times J$--diagonal matrix
  $\Delta_{J\times J}(\varphi_1)$ with $\varphi_1$ in every diagonal
  entry, and hence $g_{22}\,\varphi_1^{(J)}$ is given by the matrix
  \begin{displaymath}
    (r_{ij})_{i,j \in J} \ \Delta_{J\times
    J}(\varphi_1) = (r_{ij}\varphi_1)_{i,j \in J}.
  \end{displaymath}
  By assumption \mbox{$g_{22}\,\varphi_1^{(J)}=\varphi_1^{(J)}$}, and
  consequently we have an equality of \mbox{$J\times J$}--matrices:
  \begin{displaymath}
    (r_{ij}\varphi_1)_{i,j \in J} = \Delta_{J\times
    J}(\varphi_1).
  \end{displaymath}
  It follows that:
  \begin{displaymath}
    r_{jj}\varphi_1 = \varphi_1 \quad \text{and} \quad r_{ij}\varphi_1 = 0
    \, , \ i \neq j.
  \end{displaymath}
  Now writing $r_{ij}=a_{ij}+b_{ij}\xi$ with $a_{ij}, b_{ij} \in k$
  and applying the maps above to $1 \in k$ we get
  \begin{displaymath}
    (a_{jj}+b_{jj}\xi)\xi = \xi \quad \text{and} \quad
    (a_{ij}+b_{ij}\xi)\xi = 0 \, , \ i \neq j.
  \end{displaymath}
  We see that $a_{jj}=1$ and $a_{ij}=0$ for $i\neq j$, that is,
  \begin{displaymath}
    r_{jj}=1+b_{jj}\xi \quad \text{and} \quad r_{ij} = b_{ij}\xi
    \, , \ i \neq j.
  \end{displaymath}
  With this information at hand we can see that \mbox{$g_{22} =
    (r_{ij})_{i,j \in J}$} is invertible, in fact, $g_{22}$ is its own
  inverse. Let us simply calculate the $(i,j)$'th entry, $q_{ij}$, in
  the product matrix $g_{22}\,g_{22}$: Using that \mbox{$\xi^2=0$} and
  that the field $k \cong \mathbb{F}_2$ has characteristic $2$ it
  follows that:
  \begin{align*}
    q_{ij} &= \sum_{\nu \in J} r_{i\nu}r_{\nu j} \\
           &= \left\{ 
   \begin{array}{ccl}
     (1+b_{jj}\xi)^2 & \textnormal{for} & i=j \\
     (1+b_{ii}\xi)b_{ij}\xi + b_{ij}\xi(1+b_{jj}\xi) &
     \textnormal{for} & i \neq j   
   \end{array}
    \right. \\
           &= \left\{ 
   \begin{array}{ccl}
     1+2b_{jj}\xi & \textnormal{for} & i=j \\
     2b_{ij}\xi & \textnormal{for} & i \neq j  
   \end{array}
    \right. \\
           &= \left\{ 
   \begin{array}{ccl}
     1 & \textnormal{for} & i=j \\
     0 & \textnormal{for} & i \neq j  
   \end{array}
    \right.
  \end{align*}
  as desired.
\end{proof}

\section{Relative resolutions and Schanuel maps} \label{sec:RvsS}

\noindent
In this section we study how the resolution condition and the Schanuel
condition of Definition \eqref{dfn:conditions} are related. In
general, neither of these two conditions imply the other, however, in
Theorems \eqref{thm:StoR} and \eqref{thm:RtoS} we give necessary and
sufficient conditions for this phenomenon to happen.

\begin{dfn}
  We say that $\F$ is \emph{weakly closed under direct summands} if
  for any $F \in \F$ and any direct summand $M$ in $F$ with $F/M \in
  \F$, the module $M$ belongs to $\F$.
\end{dfn}

\begin{exa} 
  There exist precovering classes which are not weakly closed under
  direct summands:

  \medskip
  \noindent
  The precovering class $\F$ from Example \eqref{exa:nodirectsummands}
  is not closed under direct summands. As $\F$ is closed under
  set-indexed coproducts, it follows from Proposition \eqref{prp:weak}
  below that $\F$ is not weakly closed under direct summands either.
\end{exa}

\begin{prp} 
  \label{prp:weak} 
  A precovering class $\F$ is closed under direct summands if and only
  if $\,\F$ is weakly closed under direct summands and closed under
  set-indexed (respectively, countable) coproducts in
  $\mathsf{Mod}\,R$.
\end{prp}

\begin{proof}
  ``If'': Let $M$ be a direct summand of \mbox{$F \in \F$}, that is,
  there exists some module $M'$ with \mbox{$F=M \oplus M'$}. Using
  Eilenberg's swindle we consider $F^{(\mathbb{N})}$ and note that
  \begin{displaymath}
    \tag{\text{$*$}}
    M \oplus F^{(\mathbb{N})} \cong F^{(\mathbb{N})}.
  \end{displaymath}
  As $\F$ is closed under countable coproducts,
  \mbox{$F^{(\mathbb{N})} \in \F$}, and then $(*)$ implies that $M \in
  \F$ since $\F$ is weakly closed under direct summands.

  \medskip
  \noindent
  ``Only if'': If $\F$ is closed under direct summands then obviously
  $\F$ is also weakly closed under direct summands. Since $\F$ is
  precovering and closed under direct summands, the argument in
  \cite[proof of thm.\,5.4.1, $(2)\!\!\Rightarrow\!\!(1)$]{EnochsRHA}
  shows that $\F$ is closed under set-indexed coproducts.
\end{proof}

\noindent
The reason we are interested in classes which are weakly closed under
direct summands is because of the next result.

\begin{thm} \label{thm:StoR}
  A precovering class $\F$ satisfies:
  \begin{displaymath}
     \tag{\text{$\natural$}}
    (S_{M,n}) \Longrightarrow (R_{M,n}) \text{ for all modules
      $M$ and all integers $n \geqslant 0$}
  \end{displaymath}
  if and only if $\,\F$ is weakly closed under direct summands.
\end{thm}

\begin{proof}
  ``Only if'': Under the assumption of $(\natural)$ we must prove that
  $\F$ is weakly closed under direct summands.  Thus, let $M$ be a
  direct summand of a module $F$ where $F,F/M \in \F$. As
  \begin{displaymath}
    M \oplus F/M \cong 0 \oplus F
  \end{displaymath}
  we see that $M$ is $\F$--equivalent to $0$, that is,
  $\S_\F^0([M])=[M]=[0]$. Now the assumption $(\natural)$ implies the
  existence of an $\F$--resolution of $M$ of length zero,
  \begin{displaymath}
    \tag{\text{$*$}}
    0 \longrightarrow F_0 \stackrel{\partial_0}{\longrightarrow} M
    \longrightarrow 0.  
  \end{displaymath}
  As in the end of the proof of Lemma \eqref{lem:EtoR} we see that
  $\partial_0$ is an isomorphism, and hence $M \cong F_0 \in \F$ as
  desired.

  \medskip
  \noindent
  ``If'': Now assume that $\F$ is weakly closed under direct summands.
  We will prove $(\natural)$ by induction on $n \geqslant 0$.

  \medskip
  \noindent
  We begin with the case $n=0$: Suppose that $\S_\F^0([M])=[M]=[0]$.
  By definition there exist $F',F \in \F$ with $M \oplus F' \cong 0
  \oplus F \cong F$, and since $\F$ is weakly closed under direct
  summands it follows that $F_0:=M \in \F$.  Thus $M$ admits an
  $\F$--resolution of $M$ of length zero:
  \begin{displaymath}
    0 \longrightarrow F_0 \stackrel{\mathrm{id}_M}{\longrightarrow} M
    \longrightarrow 0.  
  \end{displaymath}
  \noindent
  Next we consider the case \mbox{$n>0$}: Suppose that
  \mbox{$\S_\F^n([M])=[0]$}, and take an $\F$--precover $\partial
  \colon F_0 \longrightarrow M$. By definition,
  \begin{displaymath}
    \S_\F^{n-1}([\Ker \partial]) = \S_\F^{n-1}\S_\F([M]) =
    \S_\F^n([M]) = [0], 
  \end{displaymath}
  so the induction hypothesis implies the existence of an
  $\F$--resolution of $\Ker \partial$ of length $n-1$, say,
  \begin{displaymath}
    \tag{\text{$\flat$}}
    0 \longrightarrow F_n \longrightarrow \cdots \longrightarrow F_1
    \longrightarrow \Ker \partial \longrightarrow 0.
  \end{displaymath}
  Pasting together $(\flat)$ with the complex
  \begin{displaymath}
    0 \longrightarrow \Ker \partial \longrightarrow F_0
    \stackrel{\partial}{\longrightarrow} M \longrightarrow 0 
  \end{displaymath}
  gives an $\F$--resolution of $M$ of length $n$, as desired.
\end{proof}

\begin{dfn}
  A (precovering) class $\F$ is said to be \emph{separating} if for
  every module $M \neq 0$ there exists a non-zero homomorphism $F
  \longrightarrow M$ with $F \in \F$.
\end{dfn}

\begin{lem} \label{lem:mono} 
  For a precovering class $\F$ the following hold:
  \begin{prt}
  \item If every mono $\F$--precover is an isomorphism then $\F$ is
    separating.
  \item If $\,\F$ is separating and $\partial \colon A \longrightarrow
    B$ is a homomorphism such that $\Hom_R(F,\partial)$ is mono for
    all $F \in \F$, then $\partial$ is mono.
  \end{prt}
\end{lem} 

\begin{proof}
  \proofof{(a)} Assume that every mono $\F$--precover is an
  isomorphism, and let $M$ be a module with $\Hom_R(F,M)=0$ for all $F
  \in \F$. Thus the map $0 \longrightarrow M$ is a mono
  $\F$--precover, and hence an isomorphism by assumption, that is,
  $M=0$.

  \medskip
  \noindent
  \proofof{(b)} Applying the left exact functor $\Hom_R(F,-)$, for any
  $F \in \F$, to the exact sequence,
  \begin{displaymath}
    0 \longrightarrow \Ker \partial \longrightarrow A
    \stackrel{\partial}{\longrightarrow} B  
  \end{displaymath}
  and using that $\Hom_R(F,\partial)$ is mono, we get that
  \mbox{$\Hom_R(F,\Ker \partial)=0$}. As $\F$ is separating it follows
  that $\Ker \partial=0$, that is, $\partial$ is mono.
\end{proof}

\begin{exa}
  We give two examples of precovering classes $\F$ for which there exist
  mono $\F$--precovers which are not isomorphisms:

  \medskip
  \noindent
  (a) Let $R$ be a commutative noetherian ring which is not artinian.
  As $R$ is noetherian the class $\F=\mathsf{Inj}\,R$ of injective
  $R$--modules is precovering by \cite[thm.\,5.4.1]{EnochsRHA}.
  However, as $R$ is not artinian, $\F$ is not separating by
  \cite[cor.\,2.4.11]{Xu}, and hence Lemma \eqref{lem:mono}(a) implies
  that there must exists mono $\F$--precovers which are not
  isomorphisms.

  \medskip
  \noindent
  (b) Let $R$ be a commutative integral domain, and consider for any
  module $M$ its \emph{torsion submodule},
  \begin{displaymath}
    M_T = \big\{x \in M \,|\, \text{$rx=0$ for some $r \in
      R\setminus\{0\}$} \big\}. 
  \end{displaymath} 
  A module $M$ is called \emph{torsion} if $M_T=M$, and of course the
  torsion submodule of any module is torsion.

  \medskip
  \noindent
  The torsion modules constitutes a precovering class, in fact, given
  a module $M$ it is not hard to see that the inclusion \mbox{$M_T
    \longrightarrow M$} is a torsion precover of $M$. In particular,
  \mbox{$0=R_T \longrightarrow R$} is a mono torsion precover of $R$,
  but it is not an isomorphism.
\end{exa}




\noindent
The following result shows why we are interested in precovering
classes for which every mono precover is an isomorphism.

\begin{thm} \label{thm:RtoS}
  A precovering class $\F$ satisfies: 
  \begin{displaymath}
    \tag{\text{$\flat$}}
    (R_{M,n}) \Longrightarrow (S_{M,n}) \text{ for all modules
      $M$ and all integers $n \geqslant 0$}
  \end{displaymath}
  if and only if every mono $\F$--precover is an isomorphism.
\end{thm}

\begin{proof}
  ``Only if'': Assume $(\flat)$.  We must prove that every mono
  $\F$--precover is an isomorphism.  Any mono $\F$--precover $\varphi
  \colon F_0 \longrightarrow M$ gives an $\F$--re\-so\-lu\-tion of
  $M$ of length zero:
  \begin{displaymath}
    0 \longrightarrow F_0 \stackrel{\varphi}{\longrightarrow} M
    \longrightarrow 0,
  \end{displaymath}
  and thus our assumption ensures that \mbox{$\S^0_\F([M])=[M]=[0]$}.
  This means that $M$ is a direct summand of some \mbox{$F \in \F$}
  with a quotient $F/M \in \F$.  In particular, $M$ is a homomorphic
  image of \mbox{$F \in \F$}, and this implies that the $\F$--precover
  $\varphi$ must be epi.  Consequently, $\varphi$ is an isomorphism.



  \medskip
  \noindent
  ``If'': Conversely, assume that every mono $\F$--precover is an
  isomorphism.  We must show $(\flat)$, which we do by induction on $n
  \geqslant 0$:

  \medskip
  \noindent
  We begin with the case $n=0$. Thus, let $M$ be any module for which
  there exists an $\F$--resolution of length zero:
  \begin{displaymath}
    \tag{\text{$*$}}
    0 \longrightarrow F_0 \stackrel{\partial}{\longrightarrow} M
    \longrightarrow 0. 
  \end{displaymath}
  We must argue that $\S^0_\F([M])=[0]$. Actually, we prove something
  even stronger, namely that \mbox{$M \in \F$}.  Since $(*)$ is an
  $\F$--resolution we have exactness of
  \begin{displaymath}
    \xymatrix{
    0 \ar[r] & \Hom_R(F,F_0) \ar[rr]^-{\Hom_R(F,\partial)} & & \Hom_R(F,M)
    \ar[r] & 0,
    }
  \end{displaymath}
  that is, $\Hom_R(F,\partial)$ is an isomorphism for all \mbox{$F \in
    \F$}.  Now our assumption and Lemma \eqref{lem:mono}(a) and (b)
  gives that $\partial \colon F_0 \longrightarrow M$ is a mono
  $\F$--precover.  Another application of our assumption then gives
  that $\partial$ is an isomorphism, and thus $M \cong F_0 \in \F$.

  \medskip
  \noindent
  Next we assume that \mbox{$n>0$}. Let $M$ be a
  module which has an $\F$--re\-so\-lu\-ti\-on of length $n$,
  \begin{displaymath}
    \tag{\text{$\dagger$}}
    0 \longrightarrow F_n \longrightarrow
    \cdots \longrightarrow F_1 
    \stackrel{\partial_1}{\longrightarrow} F_0
    \stackrel{\partial_0}{\longrightarrow}  M \longrightarrow 0.
  \end{displaymath}
  We break up $(\dagger)$ into two complexes,
  \begin{eqnarray}
    &0 \longrightarrow F_n \longrightarrow
    \cdots \longrightarrow F_1 
    \stackrel{\widehat{\partial}_1}{\longrightarrow} \Ker \partial_0
    \longrightarrow 0,& \\ 
     &0 \longrightarrow \Ker \partial_0
     \stackrel{\iota}{\longrightarrow} F_0
     \stackrel{\partial_0}{\longrightarrow} M \longrightarrow 0,& 
  \end{eqnarray}
  where $\widehat{\partial}_1$ is the co-restriction of $\partial_1$
  to $\Ker \partial_0$. Once we have argued that $\widehat{\partial}_1
  \colon F_1 \longrightarrow \Ker \partial_0$ is an $\F$--precover, it
  will follow that the upper sequence $(1)$ is an $\F$--resolution of
  $\Ker \partial_0$, and hence the induction hypothesis gives that
  \mbox{$\S^{n-1}_\F([\Ker \partial_0])=[0]$}. By the lower sequence
  $(2)$ we have \mbox{$\S_\F([M])=[\Ker \partial_0]$}, and thus the
  desired conclusion follows:
  \begin{displaymath}
    \S^n_\F([M])=\S^{n-1}_\F \S_\F([M])=\S^{n-1}_\F([\Ker \partial_0])=[0].
  \end{displaymath}
  To see that $\widehat{\partial}_1 \colon F_1 \longrightarrow \Ker
  \partial_0$ is an $\F$--precover we let $\varphi \colon F
  \longrightarrow \Ker \partial_0$ be any homomorphism with $F \in
  \F$. As $\partial_0\iota\varphi=0$ and as $\Hom_R(F,(\dagger))$ is
  exact, there exists $\psi \colon F \longrightarrow F_1$ with
  $\partial_1\psi = \iota\varphi$. Since $\iota$ is mono this means
  that we have a commutative diagram
  \begin{displaymath}
    \xymatrix{
      {} & F \ar[d]^-{\varphi} \ar@{-->}[dl]_-{\psi} \\
      F_1 \ar[r]_-{\widehat{\partial}_1} & \Ker \partial_0,
    }
  \end{displaymath}
  as desired.
\end{proof}

\begin{rmk}
  The dual notion of a precover is a \emph{preenvelope}, see
  \cite[chap.\,6]{EnochsRHA}. For a \emph{preenveloping} class
  $\mathsf{G}$, the reader can imagine how to construct Ext functors,
  resolutions, and Schanuel maps relative to $\mathsf{G}$, see also
  \cite[chap.\,8]{EnochsRHA}.

  \medskip
  \noindent
  Not surprisingly, every result in this this paper has an analogue in
  this ``preenveloping context''. We leave it as an exercise for the
  interested reader to verify this claim.
\end{rmk}

\providecommand{\bysame}{\leavevmode\hbox to3em{\hrulefill}\thinspace}
\providecommand{\MR}{\relax\ifhmode\unskip\space\fi MR }
\providecommand{\MRhref}[2]{%
  \href{http://www.ams.org/mathscinet-getitem?mr=#1}{#2}
}
\providecommand{\href}[2]{#2}

\bigskip
\bigskip


\begin{thebibliography}{10}

\bibitem{CE} Henri Cartan and Samuel Eilenberg, ``Homological
  algebra'', Princeton University Press, Princeton, N. J., 1956.

\bibitem{EnochsRHA} Edgar E. Enochs and Overtoun M. G. Jenda,
  ``Relative homological algebra", de Gruyter Expositions in
  Mathematics, vol.~{30}, Walter de Gruyter \& Co., Berlin, 2000.

\bibitem{EEEOMGJLO} Enochs E. Edgar, Overtoun M. G. Jenda, and Luis
  Oyonarte, \emph{{$\lambda$} and {$\mu$}-dimensions of modules},
  Rendiconti del Seminario Matematico della Universit{\`a} di Padova
  \textbf{105} (2001), 111--123.

\bibitem{IK} Irving Kaplansky, ``Commutative rings'', revised ed., The
  University of Chicago Press, Chicago, Ill.-London, 1974.

\bibitem{Xu} Jinzhong Xu, ``Flat covers of modules'', Lecture Notes in
  Mathematics, vol.~{1634}, Springer-Verlag, Berlin, 1996.

\end{thebibliography}


\end{document}